\documentclass[11pt,a4paper]{article}
\usepackage[]{times}
\usepackage{fullpage}
\usepackage{graphicx}
\usepackage{float}
\usepackage{subfigure}
\usepackage{amsmath}
\usepackage{fancyhdr}
\usepackage[colorlinks]{hyperref}
\usepackage{xcolor}
\date{}

\newcommand{\prov}{{\sc Proof}.\hspace*{0mm} }

\newcommand{\QED}{$\rule{2mm}{2mm}$}

\newtheorem{theorem}{Theorem}[section]
\newtheorem{lemma}[theorem]{Lemma}
\newtheorem{e-proposition}[theorem]{Proposition}
\newtheorem{corollary}[theorem]{Corollary}
\newtheorem{e-definition}[theorem]{Definition\rm}

\title{Fine error bounds for approximate \\ asymmetric saddle point problems}
\author{
    Vitoriano Ruas$^{1}$\thanks{Sorbonne Universit\'e, Campus Pierre et Marie Curie, 4 place jussieu, Couloir 55-65, 4\`eme \'etage, 75005 Paris, France.}
		\\[1mm]
  {\small $^{1}$ Institut Jean Le Rond d'Alembert, CNRS UMR 7190, Sorbonne Universit\'e, Paris, France.}\\[1mm]
  {\small e-mail: {\it vitoriano.ruas@upmc.fr}}}

\begin{document}
\maketitle
\vspace{-6mm}
\hspace{7cm} 
\textbf{\small Dedicated to the memory of Roland Glowinski}\\
\hspace{15mm} \textit{Declarations of interest: none}
\begin{abstract}
The theory of mixed finite element methods for solving different types of elliptic partial differential equations in saddle point formulation is well established since many decades. This topic was mostly studied for variational formulations defined upon the same product spaces of both shape- and test-pairs of primal variable-multiplier. Whenever either these spaces or the two bilinear forms involving the multiplier are distinct, the saddle point problem is asymmetric. The three inf-sup conditions to be satisfied by the product spaces stipulated in work on the subject, in order to guarantee well-posedness, are known since long (see e.g Exercise 2.14 of Ern \& Guermond (2004) \cite{ErnGuermond}). However, the material encountered in the literature addressing the approximation of this class of problems left  room for improvement and clarifications. After making a brief review of the existing contributions to the topic that justifies such an assertion, in this paper we set up finer global error bounds for the pair primal variable-multiplier solving an asymmetric saddle point problem. Besides well-posedness, the three constants in the aforementioned inf-sup conditions are identified as all that is needed for  determining the stability constant appearing therein, whose expression is exhibited. As a complement, refined error bounds depending only on these three constants are given for both unknowns separately. 


\end{abstract}

\section{Introduction}

\hspace{6mm} According to the well known theory due to Babu\v{s}ka (1971) \cite{Babuska1} and Brezzi (1974) \cite{Brezzi}, a linear variational problem (LVP) is well-posed only if the underlying continuous bilinear form is weakly coercive. In the finite dimensional case this reduces to proving that the vector space of shape-elements, as related to the one of test-elements in this bilinear form, satisfies only one \textit{inf-sup} condition. This fact was highlighted in Cuminato \& Ruas (2015) \cite{COAM}, where error bounds for  approximate LVPs were also refined.\\
\indent The particular case of saddle point problems and their approximation was studied in Babu\v{s}ka (1973) \cite{Babuska2} and Brezzi (1974) \cite{Brezzi} for symmetric problems. Here the term symmetry refers to formulations with only one bilinear form defined upon a product space for the pair primal variable-multiplier and also with the same product spaces for both shape- and test-pairs. In the asymmetric case the two product spaces are different and/or the variational formulations incorporate two distinct bilinear forms having the multiplier as an argument. Actually the above studies were extended by Nicolaides (1982) \cite{Nicolaides} to the non symmetric case. That paper should be praised as a pioneering work on the approximation of asymmetric saddle point problems. Nevertheless, it definitively left room for both clarifications and improvements, as explained hereafter. \\
\indent The mathematical background for the numerical analysis of mixed finite element methods for solving saddle point problems was well reported in the literature within the last three decades (cf. Brezzi \& Fortin (1991) \cite{BrezziFortin} or Ern \& Guermond (2004) \cite{ErnGuermond} and references therein). Most of this material is devoted to symmetric saddle point formulations. As a matter of fact, the asymmetric case has mostly been addressed in connection with specific situations, for which one-off built-in error estimates are developed. The author conjectures that this might be due to the fact that not as many mathematical models of this type are encountered in practical applications. Nevertheless, we quote a few of them. For example, in one of the first papers on the subject, MacCamy (1980) \cite{MacCamy} addressed interface problems for electromagnetic fields. Bernardi et al (1988) \cite{Bernardi} studied spectral methods for the Stokes problem using standard and weighted spaces of square integrable functions. Ciarlet Jr et al. (2003)  \cite{CiarletJr} considered a mixed formulation for the reaction-convection-diffusion equation involving a fourth bilinear form in addition to the three more common ones. 
As far as a general approximation theory for this type of problem is concerned, to the best of the author's knowledge, the aforementioned article by Nicolaides (1982) \cite{Nicolaides} is usually cited. However, in that work the verification of three additional conditions is advocated, although they are superfluous in the realistic case where the approximation spaces are finite-dimensional. Moreover, in the error bounds given therein some inf-sup constants are unnecessary, since one can do without them owing to properties of adjoint operators.\\
Actually, in this article we clarify such issues by specifying which conditions must be really checked in order to allow for more easily exploitable error bounds. Thanks to some tools of Functional Analysis (cf. Yosida (1980) \cite{Yosida}) among other ingredients, we show that there are only three and nothing else, namely, the inf-sup conditions previously identified by other authors (cf. Exercise 2.14 of Ern \& Guermond (2004) \cite{ErnGuermond}) as necessary for well-posedness, restricted here to finite-dimensional approximation spaces. The knowledge of the three underlying constants leads to optimal error estimates, in case it is proven that they are independent of the discretization level, which is known to be a rather difficult task (see e.g. Bertrand \& Ruas (2023)\cite{FBVR}).  \\ 
\indent In short, by placing our analysis in the general framework given in Cuminato \& Ruas (2015) \cite{COAM} for the approximation of LVPs, we derive proper error bounds for asymmetric saddle point problems with fine stability constants.\\ 
We should also emphasize that our error bounds are asymmetric counterparts of those established by Dupire (1985) \cite{Dupire} for symmetric saddle point problems, that we recall in Section 5. In complement to this, we further exhibit therein unified error bounds for the asymmetric case, resulting from studies also carried out in Dupire (1985) \cite{Dupire}, which give rise to best possible error estimates for the primal variable and the multiplier separately. \footnote{ In doing so, our aim is to further publish in English language for the first time in a means widely distributed with due consent (Dupire (2013) \cite{DupirePerso}), some results available only in Dupire's doctoral thesis written in Portuguese language.} Incidentally, we observe that the latter results could also be viewed as improvements of those given in Nicolaides (1982) \cite{Nicolaides} with the same purpose.\\
\indent Finally we note that, once our extension of Dupire's error bounds to the asymmetric case were made available, they could be applied to the study of a variant of Raviart-Thomas mixed method to solve elliptic PDE's in smooth domains, using different spaces of trial- and test-flux fields (cf. Bertrand \& Ruas (2023) \cite{FBVR}). This was actually an additional motivation for this work.\\

An outline of the article is as follows. Section
 2 deals with some abstract results of Functional Analysis to be used in the sequel. 
In Section 3 the analysis of the well posedness and the stability of approximate asymmetric saddle point problems is carried out. Section 4 is devoted to the global error bounds resulting from the previous study. In Section 5 we give some complementary error bounds. The main contributions of the present work are summarized in Section 6.

\section{Preparatory material} 

In this section we use the following notations. For a given Hilbert space ${\mathcal W}$ the underlying inner product is denoted by $(\cdot,\cdot)_{\mathcal W}$, the corresponding norm is denoted by $\| \cdot \|_{\mathcal W}$ and its null vector is represented by $0_{\mathcal W}$. The standard norm of a given continuous linear operator ${\mathcal A}$ between two Hilbert spaces is denoted by $\| {\mathcal A} \|$.\\
To begin with we prove
\begin{lemma} 
\label{conorm}
Let ${\mathcal V}$ and ${\mathcal U}$ be two Hilbert spaces. If ${\mathcal J}$ is a bounded linear operator from ${\mathcal V}$ onto ${\mathcal U}$ which is an isomorphism between both spaces, then we have:
$$ \displaystyle \inf_{v \in {\mathcal V} \setminus \{0_{\mathcal V}\}} \frac{\| {\mathcal J}(v) \|_{\mathcal U}}{\| v \|_{\mathcal V}} = \| {\mathcal J}^{-1} \|^{-1}. $$
\end{lemma}

\prov First we recall that in the case at hand the inverse ${\mathcal J}^{-1}$ of ${\mathcal J}$ is a well defined linear operator, which is an isomorphism between ${\mathcal V}$ and ${\mathcal U}$. Moreover from the Banach Inverse Theorem (see e.g. Yosida (1980) \cite{Yosida}), it is also a bounded operator. \\
Now, since for every $v \in {\mathcal V}$ we can define $u \in {\mathcal U}$ by $u = {\mathcal J}(v)$, we may write
\begin{equation}
\label{lowerbound}
\begin{array}{l}
\displaystyle \inf_{v \in {\mathcal V} \setminus \{0_{\mathcal V}\}} \frac{\| {\mathcal J}(v) \|_{\mathcal U}}{\| v \|_{\mathcal V}} \geq  
\displaystyle \inf_{u \in {\mathcal U} \setminus \{0_{\mathcal U}\}} \frac{\| u \|_{\mathcal U}}{\| {\mathcal J}^{-1}(u) \|_{\mathcal V}} \\
= \displaystyle \left[ \sup_{u \in {\mathcal U} \setminus \{0_{\mathcal U}\}} \frac{\| {\mathcal J}^{-1}(u) \|_{\mathcal V}}{\| u \|_{\mathcal U}} \right]^{-1} = \| {\mathcal J}^{-1} \|^{-1}. 
\end{array}
\end{equation}
Actually, on the first row of \eqref{lowerbound} the equality applies, since, conversely, for every $u \in {\mathcal U}$ $\exists v \in {\mathcal V}$ given by 
$v = {\mathcal J}^{-1}(u)$. \QED \\

Let now ${\mathcal K}$ be a bounded linear operator from $Y$ to $X$, $Y$ and $X$ being two Hilbert spaces. \\
We first recall a classical result of Functional Analysis (see e.g. Conway (1990) \cite{Conway}), namely,\\
The adjoint operator ${\mathcal K}^*$ from $X$ to $Y$ satisfies $\| {\mathcal K}^* \| = \| {\mathcal K} \|$.\\

Next we prove
\begin{e-proposition}
\label{surjection}
The operator ${\mathcal K}$ satisfies
\begin{equation}
\label{minored}
\displaystyle \| {\mathcal K} \mu \|_X \geq \kappa \| \mu \|_Y, \; \forall \mu \in Y. 
\end{equation}
for a suitable constant $\kappa >0$, if and only if its adjoint operator ${\mathcal K}^*$ is a surjection.
\end{e-proposition}

\prov 
Using the usual notation $R({\mathcal A})$ and $Ker({\mathcal A})$ for the range and the kernel of a linear operator ${\mathcal A}$, first we establish that 
\eqref{minored} holds if and only if $R({\mathcal K})$ is closed and $Ker({\mathcal K}) = \{0_Y\}$.\\

\noindent \underline{Assume that \eqref{minored} holds:} In this case we obviously have $Ker({\mathcal K}) = \{0_Y\}$. Moreover $R({\mathcal K})$ is closed due to the following argument. Let $\{\xi_n\}_n$ be a sequence of $R({\mathcal K})$ converging to $\xi \in X$. In order to prove that $\xi \in R({\mathcal K})$ we note that, by definition,  there exists a sequence $\{ \mu_n \}_n$ in $Y$ such that 
${\mathcal K} \mu_n = \xi_n$ for all $n$. Since $\{\xi_n\}_n$ is a Cauchy sequence of $X$, \eqref{minored} implies that $\{ \mu_n \}_n$ is 
a Cauchy sequence of $Y$. Therefore there exists $\mu \in Y$ such that $\mu_n$ converges to $\mu$ . By the continuity of ${\mathcal K}$, $\xi_n = {\mathcal K} \mu_n$ converges to ${\mathcal K} \mu$, which is thus equal to $\xi$.\\
\underline{Assume that $R({\mathcal K})$ is closed and $Ker({\mathcal K}) = \{0_Y\}$:} In this case ${\mathcal K}$ is an isomorphism 
between $Y$ and $R({\mathcal K})$. Hence the inverse of ${\mathcal K}$ is also bounded operator and an isomorphism between $R({\mathcal K})$ and $Y$ and we may write, $\forall \mu \in Y$, $\| \mu \|_Y = \| {\mathcal K}^{-1} {\mathcal K} \mu \|_Y \leq  
\| {\mathcal K}^{-1} \| \| {\mathcal K} \mu \|_X$. It follows that \eqref{minored} holds with $\kappa = \| {\mathcal K}^{-1} \|^{-1}$.\\

\noindent \underline{Now assume that ${\mathcal K}^*$ is a surjection:} In this case $R({\mathcal K}^*) = Y$, whose orthogonal space is $\{ 0_Y \}$, which happens to be $Ker({\mathcal K})$. Moreover since $R({\mathcal K}^*)$ is closed, by The Closed Range Theorem (cf. Yosida (1980) \cite{Yosida}) $R({\mathcal K})$ is also closed. Therefore, from the first part of the proof, \eqref{minored} holds true.\\
\underline{Now assume that \eqref{minored} holds:} In this case, the first part of the proof also tells us that $R({\mathcal K})$ is closed, and hence so is $R({\mathcal K}^*)$ again by The Closed Range Theorem. It also tells us that $Ker({\mathcal K}) = \{ 0_Y \}$, and thus the orthogonal space of $Ker({\mathcal K})$ is $Y$. Since the latter space is nothing but the closure of $R({\mathcal K}^*)$, the fact that it is closed implies that it is equal to $Y$, and the result follows. \QED \\

\begin{corollary}
\label{corollary1}
The largest possible constant $\kappa$ for which \eqref{minored} holds is $\| {\mathcal K}^{-1} \|^{-1}$, ${\mathcal K}$ being viewed as a linear operator from $Y$ onto the complete space $W:=R({\mathcal K})$.
\end{corollary}

\prov We saw that ${\mathcal K}$ is a continuous isomorphism between the Hilbert spaces $Y$ and $W$. In view of this the result is a 
direct consequence of Lemma \ref{conorm}. \QED 

\begin{corollary}
\label{corollary2}
Assuming that \eqref{minored} holds, let $Z \subset X$ be the orthogonal complement of $Ker({\mathcal K}^*)$. Then the following condition holds:
\begin{equation}
\label{minored*}
\| {\mathcal K}^* \zeta \|_Y \geq \kappa \| \zeta \|_X, \; \forall \zeta \in Z.
\end{equation}
\end{corollary}

\prov Clearly enough ${\mathcal K}^*$ is a continuous isomorphism between $Z$ and $Y$ since from Proposition \ref{surjection} it is a surjection. Therefore From Corollary \ref{corollary1} we have 
\begin{equation}
\label{lowerbound*}
\| {\mathcal K}^* \zeta \|_Y \geq  \| {\mathcal K^*}^{-1} \|^{-1} \| \zeta \|_X, \; \forall \zeta \in Z.
\end{equation}
On the other hand it is well known that ${\mathcal K^*}^{-1} = {\mathcal K^{-1}}^{*}$ and thus $\| {\mathcal K^*}^{-1} \| = 
\| {\mathcal K^{-1}}^{*} \| = \|{\mathcal K^{-1}} \|$. Since from Corollary \ref{corollary1} $\kappa \leq \|{\mathcal K^{-1}} \|^{-1}$ the 
result follows. \QED \\

Before proceeding to the study of asymmetric mixed formulations, we recall a lemma due to Dupire (1985) \cite{Dupire}, whose proof we give in full here, in order to make the arguments that follow self-contained. 
\begin{lemma}
\label{LemmaII-5-1}
(cf. Dupire (1985) \cite{Dupire}) \\ Let $(a_1,a_2,a_3) \in (\Re^+)^3$ and $a_4 := \max[a_1x_1-a_2x_2,a_3x_2]$. Then $\forall (x_1;x_2) \in (\Re^+)^2$ we have $a_4 \geq a_5 \sqrt{x_1^2+x_2^2}$, where $a_5:= a_1 a_3/\sqrt{a_1^2+(a_2 + a_3)^2}$.
\end{lemma}
\prov
We have $x_2 \leq a_4/a_3$ and $x_1 \leq (a_4+a_2x_2)/a_1 \leq a_4 (a_2+a_3)/a_1 a_3$. Thus $x_1^2+x_2^2 \leq a_4^2 \{[(a_2+a_3)/(a_1 a_3)]^2 + 1/a_3^2\}$ and the result follows. \QED 
 
\section{Weak coercivity of asymmetric saddle-point bilinear forms}

\hspace{6mm} This paper's abstract setting deals with bilinear forms defined upon two pairs of product spaces, namely, $W:=P \times U$ and $Z:=Q \times V$, where $P$, $Q$, $U$ and $V$ are four finite-dimensional spaces fulfilling $dim P=dim Q$ and $dim U =dim V$ equipped with their respective norms and inner products denoted as indicated in the previous section. The product spaces in turn are equipped with the norms $|| (p;u) ||_W :=  [\| p \|_P^2 + \| u \|_U^2]^{1/2}$ $\forall (p;u) \in W$ and $|| (q;v) ||_Z :=  [\| q \|_Q^2 + \| v \|_V^2]^{1/2}$ $\forall (q;v) \in Z$. \\
More precisely, we consider bilinear forms given by 
\begin{equation}
\label{bilinearformc}
c((p;u),(q;v)):= a(p,q) + b(u,q) + d(p,v) \; \forall (p;u) \in W \mbox{ and } \forall (q;v) \in Z,
\end{equation}
where $a$, $b$ and $d$ are continuous bilinear forms defined upon $P \times Q$, $P \times V$ and $U \times Q$, respectively. \\
Henceforth we denote by $\|a\|$, $\|b\|$ and $\| d \|$ three constants viewed as the norms of $a$, $b$ and $d$ fulfilling $a(p,q) \leq \|a\| \|p\|_P \| q \|_Q \; \forall (p;q) \in P \times Q$, 
$b(u,q) \leq \|b\| \|u\|_U \|q\|_Q \; \forall (u;q) \in U \times Q$ and $d(p,v) \leq \|d\| \|p\|_P \|v\|_V$ $\forall (p;v) \in P \times V$. It is easy to see that $c$ is also continuous over $W \times Z$ with $\|c\|:=\max[\|a\|,\|b\|,\|d\|]$. \\  
Let $B: Q \rightarrow U$ and $D:P \rightarrow V$ be the continuous linear operators defined by
\begin{equation}
\label{operators}
\begin{array}{ll} 
(u,Bq)_U= b(u,q) \; \forall (u;q) \in U \times Q; & (Dp,v)_V= d(p,v) \; \forall (p;v) \in P \times V.
\end{array}
\end{equation}
Notice that the norm $\| B \|$ is bounded above by $\| b\|$. Indeed, clearly enough, $\forall q \in Q$
$$\| B q \|_U = \displaystyle \sup_{u \in U \setminus \{0_U\}} \frac{(u,Bq)_U}{\| u \|_U} 
= \displaystyle \sup_{u \in U \setminus \{0_U\}} \frac{b(u,q)}{\| u \|_U} \leq \| b \| \| q \|_Q.$$ 
 $B^{*}: U \rightarrow Q$ is defined by 
Similarly, the norm of the operator $D$ is bounded above by $\|d\|$. Moreover the norm of the adjoint operator $D^{*}: V \rightarrow P$ defined by 
$(p,D^{*}v)_P = d(p,v) \; \forall p \in P \mbox{ and } \forall v \in V$ 
is also bounded above by $\| d \|$.\\ 

Next we endeavor to exhibit conditions on $a$, $b$ and $d$ such that $c$ is continuous and weakly coercive. \\
First of all, the continuity of $c$ means that there exists a constant $\|c\|$, viewed as the norm of $c$, such that
\begin{equation}
\label{cc}
c((p;u),(q;v)) \leq \|c\| \|(p;u)\|_W \|(q;v)\|_Z \; \forall (p;u) \in W \mbox{ and } \forall (q;v) \in Z.
\end{equation}
On the other hand, in the context of finite-dimensional spaces (see e.g. Cuminato \& Ruas (2015) \cite{COAM}), the weak coercivity of $c$ is ensured solely if there exists a constant $\gamma>0$ such that 
\begin{equation}
\label{wcc}
\forall (p;u) \in W \displaystyle \sup_{(q;v) \in Z \setminus \{0_Z\}} \frac{c((p;u),(q;v))}{\| (q;v) \|_Z} \geq \gamma \| (p;u) \|_W.
\end{equation}  
Next we prove the following fundamental result.
\begin{theorem}
\label{theo}
Let $R:=Ker(D)$, $S:=Ker(B)$ and $M$ and $N$ be the orthogonal complements of $R$ and $S$ in $P$ and $Q$, respectively. We assume that $dim R = dim S$, and also that the forms $a$, $b$ and $d$ satisfy the following \textit{inf-sup} conditions:\\
There exist real numbers $\alpha >0$, $\beta > 0$ and $\delta >0$ such that
\begin{equation}
\label{infsups}
\left\{
\begin{array}{l}
\forall r \in R \displaystyle \sup_{s \in S \mbox{ s.t. } \| s \|_Q=1} a(r,s) \geq \alpha \| r \|_P \\
\forall u \in U \displaystyle \sup_{q \in Q \mbox{ s.t. } \| q \|_Q=1} b(u,q) \geq \beta  \| u \|_U \\
\forall v \in V \displaystyle \sup_{p \in P \mbox{ s.t. } \| p \|_P=1} d(p,v) \geq \delta \| v \|_V.
\end{array}
\right.
\end{equation}
Then there exists a constant $\gamma>0$ such that \eqref{wcc} holds.
\end{theorem}
\prov
The second \textit{inf-sup} condition in \eqref{infsups} implies that  
\begin{equation}
\label{existenceB}
\forall u \in U \; \exists q \in Q \mbox{ with } \| q \|_Q = 1 \mbox{ s.t. } 
b(u,q) \geq \beta  \| u \|_U.
\end{equation}
Analogously, the third \textit{inf-sup} condition in \eqref{infsups} implies that 
\begin{equation}
\label{lowerboundD}
\| D^{*} v \|_P := \displaystyle \sup_{p \in P \mbox{ s.t. } \| p \|_P=1} d(p,v) \geq \delta \| v\|_V \; \forall v \in V.
\end{equation}
Now, setting ${\mathcal K}= D^{*}$, $Y = U$ and $X=M$, since $[D^*]^* = D$ Corollary \ref{corollary2} allows us to assert that,
\begin{equation}
\label{lowerboundortho}
\| D m \|_V := \displaystyle \sup_{v \in V \mbox{ s.t. } \| v \|_V=1} d(m,v) \geq \delta \| m \|_P \; \forall m \in M.
\end{equation}
It is noticeable that \eqref{lowerboundortho} is equivalent to the condition 
\begin{equation}
\label{existenceortho}
\forall m \in M \; \exists v \in V \mbox{ with } \| v \|_V=1 \mbox{ s.t. } 
d(m,v) \geq \delta \| m \|_P.
\end{equation} 
Now we proceed as follows: \\
Given a pair $(p; u) \in W$, we associate with it three pairs $(q_j;v_j) \in Z$, fulfilling $\| (q_j; v_j) \|_Z = 1$ for $j=1,2,3$ defined in accordance to \eqref{existenceB}, \eqref{existenceortho} and the first \textit{inf-sup} condition in \eqref{infsups}. More precisely, we first define the following fields: 
\begin{itemize}
\item 
$r \in R$ is the orthogonal projection of $p$ onto $R$;
\item
$m \in M$ is the orthogonal projection of $p$ onto $M$.
\end{itemize}
Then we set
\begin{enumerate}
\item 
$v_1 = 0_{V}$ and $q_1$ fulfills $b(u,q_1) \geq \beta \| u \|_U$ with $\| q_1 \|_Q =1$. \\
It is easy to see that $c((p;u),(q_1;v_1)) = a(p,q_1) + b(u,q_1) \geq  \beta \| u \|_U - \|a\| \| p \|_P$;
\item
$q_2 = 0_{Q}$ and $v_2 \in V$ is the function associated with $m$ fulfilling \eqref{existenceortho}, that is, $d(m,v_2) \geq \delta \| m \|_P$with $\| v_2 \|_V = 1$. \\
Since $p=r+m$, it is also clear that $c((p;u),(q_2;v_2)) = d(m,v_2) \geq \delta \| m \|_P$; 
\item 
$v_3 = 0_V$, and $q_3 \in S$ fulfills $a(r,q_3) \geq \alpha \| r \|_{P}$ with $\| q_3 \|_Q =1$.\\
By a straightforward calculation we have $c((p;u),(q_3;v_3)) = a(r,q_3) + a(m,q_3) \geq \alpha \| r \|_P - \| a \| \| m \|_P$.
\end{enumerate}
Taking into account the definitions and properties listed above, we have
\begin{equation}
\label{lowerboundc}
\displaystyle \sup_{(q;v) \in Z \mbox{ s.t. } || (q;v) ||_Z = 1}
c((p; u),(q;v)) \geq \displaystyle \max_{1 \leq i \leq 3} c((p;u),(q_i;v_i)) \geq \Phi(\| r \|_P, 
\| m \|_P,\| u \|_U),
\end{equation}
where $\Phi(x,y,z) := \max[ \beta z - \| a \| t, \delta y, \alpha x - \| a \| y)$ with $t= \sqrt{x^2+y^2}$, $x,y,z$ being non negative real numbers.\\
Using Lemma \ref{LemmaII-5-1} we can write, $max[\delta y, \alpha x - \| a \| y \geq \nu t$ with  
\begin{equation}
\label{deltac}
\nu = \alpha \delta/\sqrt{(\| a \|+\delta)^2  + \alpha^2}.
\end{equation}
Hence we can write $\Phi(x,y,z) \geq \max[\beta z - \| a \|t, \nu t]$, which yields
\begin{equation}
\label{gamma} 
\Phi(x,y,z) \geq \gamma \sqrt{x^2+y^2+z^2} \mbox{ with } \gamma:= \beta \nu/\sqrt{(\nu +\| a \|)^2+ \beta^2},
\end{equation}
by applying again Lemma \ref{LemmaII-5-1}. \\
Plugging \eqref{gamma} into \eqref{lowerboundc}, we obtain 
\begin{equation}
\label{conormc}
\displaystyle \sup_{(q;v) \in Z \mbox{ s.t. } \| (q;v) \|_Z = 1}
c((p;u),(q;v)) \geq \gamma \sqrt{\|r\|_P^2 + \| m \|_P^2 + \| u \|_U^2}.
\end{equation}
It follows from \eqref{conormc} and the Pythagorean Theorem that \eqref{wcc} holds true with $\gamma$ defined by \eqref{gamma}-\eqref{deltac}.
\QED \\

\section{Error bounds}

\hspace{6mm} It is easy to figure out that if one of the inf-sup conditions in \eqref{infsups} does not hold, the bilinear form $c$ cannot be weakly coercive. Conversely, as we are dealing with finite-dimensional spaces, nothing else is needed for $c$ to be weakly coercive (cf. Cuminato \& Ruas (2015) \cite{COAM}). 
Hence, taking into account Theorem \ref{theo}, if the uniform stability of a mixed finite element approximate problem based on an asymmetric saddle point formulation is to be established, the following issues must be overcome: \\
\indent First of all it is necessary to prove the validity of the inf-sup conditions in \eqref{infsups} for the three bilinear forms $a$, $b$ and $d$ involved in the larger bilinear form $c$ defined upon two different finite-element product spaces. This will guarantee the well-posedness of the approximate problem. Once these results are obtained, uniform stability will be a consequence of the fact that the constants $\alpha$, $\beta$ and $\delta$ appearing in \eqref{infsups} do not depend on discretization parameters such as mesh sizes.\\ 
Keeping this in mind, we next address the error bound for the solution of a saddle point problem, viewed as an approximation of another pair of unknowns $(\tilde{p};\tilde{u})$ in the following framework:\\
Without loss of generality, we consider that $(\tilde{p};\tilde{u})$ belongs to a space of infinite dimension $\tilde{W}:=\tilde{P} \times \tilde{U}$, where $\tilde{P}$ (resp. $\tilde{U}$) is equipped with an inner product and contains $P$ (resp. $U$). $F$ being a bounded linear functional on another infinite-dimensional product space $\tilde{Z}:= \tilde{Q} \times \tilde{V}$, where $\tilde{Q}$ (resp. $\tilde{V}$) is equipped with an inner product and contains $Q$ (resp. $V$), the approximation of $(\tilde{p};\tilde{u})$ is the solution $(p;u) \in W$ of a problem of the form
\begin{equation}
\label{approx}
c((p;u),(q;v))=F((q;v)) \; \forall (q;v) \in Z.
\end{equation}
\indent Assuming for simplicity that a conforming method is being used, $(\tilde{p};\tilde{u})$ in turn is the solution of a saddle point problem of the same nature in the space $\tilde{W}$ with test pairs $(\tilde{q};\tilde{v}) \in \tilde{Z}$. This means that   
 $c$ is also defined and continuous on $\tilde{W} \times \tilde{Z}$ with the same norm $\|c\|$ and that the inner product and norm $(\cdot,\cdot)_{\tilde{W}}$ and $\|\cdot\|_{\tilde{W}}$ (resp. $(\cdot,\cdot)_{\tilde{Z}}$ and $\|\cdot\|_{\tilde{Z}}$) restricted to $W$ (resp. $Z$) reduces to $(\cdot,\cdot)_W$ and $\|\cdot\|_W$ (resp. $(\cdot,\cdot)_Z$ and $\|\cdot\|_Z$). \\
\indent Resorting to the theoretical setting of Cuminato \& Ruas (2015) \cite{COAM}, it holds
\begin{equation}
\label{distance}
\| (\tilde{p};\tilde{u}) - (p;u) \|_{\tilde W} \leq  \gamma^{-1} \|c\| \displaystyle \inf_{(r;y) \in W} \| (\tilde{p};\tilde{u}) - (r;y) \|_{\tilde{W}}.
\end{equation}
\indent Inequality \eqref{distance} tells us that the existence of a lower bound for the weak coercivity constant $\gamma$ of the bilinear form $c$ independent of the dimension of approximation spaces is essentially all that is needed to attain the best possible global error estimates for the conforming method in use. According to \eqref{gamma}-\eqref{deltac}, this property is a consequence of the existence of uniform lower bounds for the constants $\alpha$, $\beta$ and $\delta$. \\
 
\section{Complementary results}

\hspace{6mm} As a complement to the above studies, we consider some possible refinements of the error bound \eqref{distance}, the proofs of which can be found in Dupire (1985) \cite{Dupire}.\\

First of all, in the symmetric case the bilinear forms $b$ and $d$ and the spaces $U \times P$ and $V \times Q$ coincide. Therefore $\delta = \beta$ and in this case the weak-coercivity
 constant $\gamma$ is expressed as follows:
\begin{equation}
\label{symmetric}
\begin{array}{l}
\gamma:= \displaystyle \frac{\alpha \beta^2}{\sqrt{(\alpha \beta + \lambda \| a \|)^2+ \lambda^2 \beta^2}} \\
\mbox{with} \\
\lambda:= \sqrt{(\|a\|+\beta)^2+ \alpha^2}
\end{array}
\end{equation}
Since $\lambda > \sqrt{\|a\|^2+\beta^2}$ if the method in use is stable, after straightforward calculations we easily come up with
\begin{equation}
\label{gammaup}
\gamma < \bar{\gamma} := \displaystyle \frac{\alpha \beta^2}{\|a\|^2+\beta^2}.
\end{equation}
It turns out that the bilinear form $c$ defined on $(V \times Q) \times (V \times Q)$ also admits $\bar{\gamma}$ as a larger weakly coercive constant. Moreover $\bar{\gamma}$ is expressed in a simpler way than $\gamma$. However, since $\alpha \leq \|a \|$, both constants scale rather similarly in terms of the key parameters $\|a\|$, $\alpha$ and $\beta$. For this reason we refrain from considering any adaptation to the asymmetric case of the above refinement of the weak-coercivity constant for the bilinear form $c$. \\

In contrast, the following error bounds for approximations of asymmetric saddle point problems can be very useful, in case the order of the approximation of the primal variable in the norm of $P$ is not the same as the one of the multiplier in the norm of $U$. As already stated in Section 1, the material that follows can also be viewed as a refinement of the results given in Nicolaides (1982) \cite{Nicolaides}. 
 \\
\indent To begin with, we recall below the definition of the (sinus of the) angle $\Gamma({\mathcal W}_1,{\mathcal W}_2)$ between two closed subspaces ${\mathcal W}_1$ and ${\mathcal W}_2$ of a normed space ${\mathcal W}$ with norm $\| \cdot \|_{\mathcal W}$ (see e.g. Kato (1966) \cite{Kato})   
\begin{equation}
\label{angle}
\Gamma({\mathcal W}_1,{\mathcal W}_2) = \displaystyle \sup_{w_1 \in {\mathcal W}_1 \mbox{ such that } \| w_1 \|_{\mathcal W} =1} \inf_{w_2 \in {\mathcal W}_2} \| w_1-w_2 \|_{\mathcal W} 
\end{equation}
It is noticeable that in case one of the spaces ${\mathcal W}_1$ and ${\mathcal W}_2$ is a subspace of the other then $\Gamma({\mathcal W}_1,{\mathcal W}_2)=0$.\\
\indent In the remainder of this section, following Dupire (1985) \cite{Dupire}, we assume that $a$, $b$ and $d$ are also continuous bilinear forms upon the product spaces considered in the previous section, namely, $\tilde{P} \times \tilde{Q}$, $\tilde{U} \times \tilde{Q}$ and $\tilde{P} \times \tilde{V}$ with the same norms $\|a\|$, $\|b\|$, $\|d\|$. \\
\indent Let $\tilde{R}$ be the kernel of operator $D:\tilde{P} \rightarrow \tilde{V}$, that is, the subspace of $\tilde{P}$ defined as follows:
$$\tilde{R} := \{\tilde{r} \; | d(\tilde{r},\tilde{v})=0 \; \forall \tilde{v} \in \tilde{V}\}.$$
Recalling the kernel $R$ of the projection onto $V$ of the same operator restricted to $P$, first we give an error bound for the approximation $p$ of the exact primal variable $\tilde{p}$, namely,
\begin{equation}
\label{errorprimal} 
\| \tilde{p} - p \|_{\tilde P} \leq  \displaystyle \frac{\|a\| \| d\|}{\alpha \delta}  \displaystyle \inf_{r \in P} \| \tilde{p} - r \|_{\tilde{P}}
+ \frac{\|b\|}{\alpha} \Gamma(R,\tilde{R}) \displaystyle \inf_{y \in U} \| \tilde{u} - y \|_{\tilde{U}}.
\end{equation}
An important consequence of \eqref{errorprimal} is the fact that the quality of the approximation of the multiplier $\tilde{u}$ in the subspace $U$ does not interfere in the quality of the approximation $p$ of the exact primal variable $\tilde{p}$, as long as the kernel space $R$ is a subspace of the kernel subspace $\tilde{R}$. This situation does occur for well known approximations methods, such as Raviart-Thomas' for simplexes (see e.g. Brezzi \& Fortin (1991) \cite{BrezziFortin}). Notice also that, in the general case, it is likely that the two kernels spaces $R$ and $\tilde{R}$ are very close to each other, so that the angle between them is small. Hence, in this case the influence of the interpolation error for the multiplier in its approximation space $U$ will diminish, even though a lower order of this error might prevail yet in the global error of $\tilde{p}$. \\ 
Regarding an error bound for the multiplier, we have
\begin{equation}
\label{errormultipli} 
\| \tilde{u} - u \|_{\tilde U} \leq  \displaystyle \frac{\|a\|^2 \| d\|}{\alpha \beta \delta}  \displaystyle \inf_{r \in P} \| \tilde{p} - r \|_{\tilde{P}} + \frac{\|b\|}{\beta} \displaystyle \left[ 1+ \frac{\| a \|}{\alpha} \Gamma(R,\tilde{R}) \right] \displaystyle \inf_{y \in U} \| \tilde{u} - y \|_{\tilde{U}}.
\end{equation}
\eqref{errormultipli} is nothing but a confirmation that the error of the multiplier is always influenced by the quality of the approximation properties of both spaces $P$ and $U$, as already pointed out by other authors (see e.g. Ern \& Guermond (2004) \cite{ErnGuermond}).\\

\section{Conclusion}

\hspace{6mm}  
The results given here are refinements or clarifications of existing ones,   aimed at providing a more comprehensive tool for the numerical analysis of the class of problems at hand. More concretely, the author can summarize as follows the main contributions of this work. 

\begin{enumerate}
\item 
First of all we stressed the fact that, if the approximate problem is posed as usual in terms of finite-dimensional subspaces, only three (inf-sup) conditions must be satisfied for it to be well-posed.
\item
Secondly we clarified that, if an inf-sup condition holds with a certain constant for any of the two bilinear forms having as arguments the multiplier and a primal variable in the orthogonal of the kernel of the associated operator, then the inf-sup condition by inverting the arguments will also hold with the same constant.
\item 
Most important, we combined the unified framework for LVPs of Cuminato \& Ruas (2015) \cite{COAM} and Theorem 3.1 to obtain error bounds with finer constants than in classical work on the topic, whose quantity is also reduced to the strictly necessary. 
\end{enumerate} 
\hspace{6mm} In short, more proper error bounds, as compared to previous studies on approximate problems of the form \eqref{approx} with a bilinear form given by \eqref{bilinearformc}, were established for conforming methods. The non conforming case among other variational crimes can be dealt with in the way advocated for instance in Cuminato \& Ruas (2015) \cite{COAM}. As far as the author can see, the constants to determine in order to apply these error bounds are close to the best possible ones. This is a relevant asset of this work, for such formulations arise more often in practice than one might think, as shown by some references listed below, including the author's joint work with Bertrand (2023) \cite{FBVR}.

\section*{Declaration}

The author declares no competing interests.

\end{document}